\documentclass[11pt]{amsart}

\usepackage{graphicx}
\usepackage{epsfig}

\usepackage{amsmath}
\usepackage{amssymb}
\usepackage{bigints}
\usepackage{tikz}
\usepackage{graphicx}
\usepackage{epsfig}

\usepackage{multicol}
\usepackage{xcolor}

\usetikzlibrary{shapes,arrows} 

\usetikzlibrary{decorations.markings}
\tikzset{
    set arrow inside/.code={\pgfqkeys{/tikz/arrow inside}{#1}},
    set arrow inside={end/.initial=>, opt/.initial=},
    /pgf/decoration/Mark/.style={
        mark/.expanded=at position #1 with
        {
            \noexpand\arrow[\pgfkeysvalueof{/tikz/arrow inside/opt}]{\pgfkeysvalueof{/tikz/arrow inside/end}}
        }
    },
    arrow inside/.style 2 args={
        set arrow inside={#1},
        postaction={
            decorate,decoration={
                markings,Mark/.list={#2}
            }
        }
    },
}

\theoremstyle{plain}
\newtheorem{theorem}                {Theorem}      [section]
\newtheorem*{theorem*}                {Theorem \ref{thm:appl}}
\newtheorem{proposition}  [theorem]  {Proposition}
\newtheorem{corollary}    [theorem]  {Corollary}

\theoremstyle{definition}

\newtheorem{remark}       [theorem]  {Remark}

\newenvironment{sproof}{\proof}{\endproof}

\DeclareMathOperator{\trace}{trace} 

\DeclareMathOperator{\Div}{div} 
 
\DeclareMathOperator{\ricci}{Ricci}

\DeclareMathOperator{\cst}{constant}

\DeclareMathOperator{\grad}{grad}

\DeclareMathOperator{\Graf}{Graf}
\numberwithin{equation}{section}

\begin{document}

\title[Global properties of biconservative surfaces in $\mathbb{R}^3$ and $\mathbb{S}^3$]
{Global properties of biconservative surfaces \\ in $\mathbb{R}^3$ and $\mathbb{S}^3$}

\author{Simona~Nistor}
\author{Cezar~Oniciuc}

\address{Faculty of Mathematics - Research Department\\ Al. I. Cuza University of Iasi\\
Bd. Carol I, 11 \\ 700506 Iasi, Romania} \email{nistor.simona@ymail.com}

\address{Faculty of Mathematics\\ Al. I. Cuza University of Iasi\\
Bd. Carol I, 11 \\ 700506 Iasi, Romania} \email{oniciucc@uaic.ro}


\thanks{The authors' work was supported by a grant of the Romanian National Authority for Scientific Research and Innovation, CNCS - UEFISCDI, project number PN-II-RU-TE-2014-4-0004.}

\subjclass[2010]{Primary 53A10; Secondary 53C40, 53C42}

\keywords{Biconservative surfaces, complete surfaces, mean curvature function, real space forms, minimal surfaces}

\begin{abstract}
We survey some recent results on biconservative surfaces in $3$-di\-mension\-al space forms $N^3(c)$ with a special emphasis on the $c=0$ and $c=1$ cases. We study the local and global properties of such surfaces, from extrinsic and intrinsic point of view. We obtain all non-$CMC$ complete biconservative surfaces in $\mathbb{R}^3$ and $\mathbb{S}^3$.
\end{abstract}

\maketitle
\section{Introduction}

The study of submanifolds with constant mean curvature, i.e., $CMC$ submanifolds, and, in particular, that of $CMC$ surfaces in $3$-dimensional spaces, represents a very active research topic in Differential Geometry for more than $50$ years.

There are several ways to generalize these submanifolds. For example, keeping the $CMC$ hypothesis and adding other geometric hypotheses to the submanifold or, by contrast, in the particular case of hypersurfaces in space forms, studying the hypersurfaces which are ``highly non-$CMC$''.

The biconservative submanifolds seem to be an interesting generalization of $CMC$ submanifolds. Biconservative submanifolds in arbitrary manifolds (and in particular, biconservative surfaces) which are also $CMC$ have some remarkable properties (see, for example \cite{FOP,MOR16-2,N-B17,S}). $CMC$ hypersurfaces in space forms are trivially biconservative, so more interesting is the study of biconservative hypersurfaces which are non-$CMC$; recent results in non-$CMC$ biconservative hypersurfaces were obtained in \cite{FT,MOR16,N,T15,UT16}.

The biconservative submanifolds are closely related to the biharmonic submanifolds. More precisely, let us consider the \textit{bienergy functional} defined for all smooth maps between two Riemannian manifolds $\left(M^m,g\right)$ and $\left(N^n,h\right)$ and given by
$$
E_{2}(\varphi)=\frac{1}{2}\int_{M}|\tau(\varphi)|^{2}\ v_g, \qquad \varphi\in C^{\infty}(M,N),
$$
where $\tau(\varphi)$ is the tension field of $\varphi$. A critical point of $E_2$ is called a \textit{biharmonic map} and is characterized by the vanishing of the \textit{bitension field} $\tau_2(\varphi)$ (see \cite{GYJ}).

A Riemannian immersion $\varphi:M^m\to \left(N^n,h\right)$ or, simply, a submanifold $M$ of $N$, is called \textit{biharmonic} if $\varphi$ is a biharmonic map.

Now, if $\varphi:M\rightarrow(N,h)$ is a fixed map, then $E_2$ can be thought as a functional defined on the set of all Riemannian metrics on $M$. This new functional's critical points are Riemannian metrics determined by the vanishing of the {\it stress-bienergy tensor} $S_2$. This tensor field satisfies
$$
\Div S_2=-\langle\tau_2(\varphi),d\varphi\rangle.
$$
If $\Div S_2=0$ for a submanifold $M$ in $N$, then $M$ is called a biconservative submanifold and it is characterized by the fact that the tangent part of its bitension field vanishes. Thus we can expect that the class of biconservative submanifolds to be much larger than the class of biharmonic submanifolds.

The paper is organized as follows. After a section where we recall some notions and general results about biconservative submanifolds, we present in \textit{Section $3$} the local, intrinsic characterization of biconservative surfaces. The local, intrinsic characterization theorem provides the necessary and sufficient conditions for an abstract surface $\left(M^2,g\right)$ to admit, locally, a biconservative embedding with positive mean curvature function $f$ and $\grad f\neq 0$ at any point.

Our main goal is to extend the \textit{local} classification results for biconservative surfaces in $N^3(c)$, with $c=0$ and $c=1$, to \textit{global} results, i.e., we ask that biconservative surfaces to be \textit{complete}, with $f>0$ everywhere and $|\grad f|>0$ on an open dense subset.

In \textit{Section $4$} we consider the global problem and construct complete biconservative surfaces in $\mathbb{R}^3$ with $f>0$ on $M$ and $\grad f\neq 0$ at any point of an open dense subset of $M$. We determine such surfaces in two ways. One way is to use the local, extrinsic characterization of biconservative surfaces in $\mathbb{R}^3$ and ``glue'' two pieces together in order to obtain a complete biconservative surface. The other way is more analytic and consists in using the local, intrinsic characterization theorem in order to obtain a biconservative immersion from $\left(\mathbb{R}^2,g_{C_0}\right)$ in $\mathbb{R}^3$ with $f>0$ on $\mathbb{R}^2$ and $|\grad f|> 0$ on an open dense subset of $\mathbb{R}^2$ (the immersion has to be unique); here, $C_0$ is a positive constant and therefore we obtain a one-parameter family of solutions. It is worth mentioning that, by a simple transformation of the metric $g_{C_0}$, $\left(\mathbb{R}^2,\sqrt{-K_{C_0}}g_{C_0}\right)$ is (intrinsically) isometric to a helicoid.

In the \textit{last section} we consider the global problem of biconservative surfaces in $\mathbb{S}^3$ with $f>0$ on $M$ and $\grad f\neq 0$ at any point of an open dense subset of $M$. As in the $\mathbb{R}^3$ case, we use the local, extrinsic classification of biconservative surfaces in $\mathbb{S}^3$, but now the ``gluing'' process is not as clear as in $\mathbb{R}^3$. Further, we change the point of view and use the local, intrinsic characterization of biconservative surfaces in $\mathbb{S}^3$. We determine the complete Riemannian surfaces  $\left(\mathbb{R}^2,g_{C_1,C^\ast_1}\right)$ which admit a biconservative immersion in $\mathbb{S}^3$ with $f>0$ everywhere and $|\grad f|>0$ on an open dense subset of $\mathbb{R}^2$ and we show that, up to isometries, there exists only a one-parameter family of such Riemannian surfaces indexed by $C_1$.

We end the paper with some figures, obtained for particular choices of the constants, which represent the non-$CMC$ complete biconservative surfaces in $\mathbb{R}^3$ and the way how these surfaces can be obtained in $\mathbb{S}^3$.

\section{Biconservative submanifolds; general properties}

Throughout this work, all manifolds, metrics, maps are assumed to be smooth, i.e. in the $C^\infty$ category, and we will often indicate the various Riemannian metrics by the same symbol $\langle,\rangle$. All surfaces are assumed to be connected and oriented.

A {\it harmonic map} $\varphi:\left(M^m,g\right)\rightarrow\left(N^n,h\right)$ between two Riemannian manifolds is a critical point of the \textit{energy functional}
$$
E:C^{\infty}(M,N)\rightarrow\mathbb{R},\quad E(\varphi)=\frac{1}{2}\int_{M}|d\varphi|^{2}\ v_g,
$$
and it is characterized by the vanishing of its {\it tension field}
$$
\tau(\varphi)=\trace_g \nabla d\varphi.
$$

The idea of the stress-energy tensor associated to a functional comes from D. Hilbert (\cite{H}). Given a functional $E$, one can associate to it a symmetric 2-covariant tensor field $S$ such that $\Div S=0$ at the critical points of $E$. When $E$ is the energy functional, P. Baird and J. Eells (\cite{BE}), and A. Sanini (\cite{S83}), defined the tensor field
$$
S=e(\varphi)g-\varphi^\ast h=\frac{1}{2}|d\varphi|^2g-\varphi^\ast h,
$$
and proved that
$$
\Div S=-\langle\tau(\varphi),d\varphi\rangle.
$$
Thus, $S$ can be chosen as the stress-energy tensor of the energy functional. It is worth mentioning that $S$ has a variational meaning. Indeed, we can fix a map $\varphi:M^m\to\left(N^n,h\right)$ and think $E$ as being defined on the set of all Riemannian metrics on $M$. The critical points of this new functional are Riemannian metrics determined by the vanishing of their stress-energy tensor $S$.

More precisely, we assume that $M$ is compact and denote
$$
\mathcal{G}=\left\{g\ : \ g \text{ is a Riemannian metric on } M\right\}.
$$
For a deformation $\left\{g_t\right\}$ of $g$ we consider $\omega=\left.\frac{d}{dt}\right|_{t=0}g_t\in T_g\mathcal{G}=C\left(\odot^2T^\ast M\right)$. We define the new functional
$$
\mathcal{F}:\mathcal{G}\to\mathbb{R}, \quad \mathcal{F}(g)=E(\varphi)
$$
and we have the following result.

\begin{theorem}[\cite{BE,S83}]
Let $\varphi:M^m\rightarrow\left(N^n,h\right)$ and assume that $M$ is compact. Then
$$
\left.\frac{d}{dt}\right|_{t=0} \mathcal{F}\left(g_t\right)=\frac{1}{2}\int_M \langle \omega,e(\varphi)g-\varphi^\ast h\rangle \ v_g.
$$
Therefore $g$ is a critical point of $\mathcal{F}$ if and only if its stress-energy tensor $S$ vanishes.
\end{theorem}

We mention here that, if $\varphi:\left(M^m,g\right)\to \left(N^n,h\right)$ is an arbitrary isometric immersion, then $\Div S=0$.

A natural generalization of harmonic maps is given by biharmonic maps. A {\it biharmonic map} $\varphi:\left(M^m,g\right)\rightarrow\left(N^n,h\right)$ between two Riemannian manifolds is a critical point of the \textit{bienergy functional}
$$
E_2:C^{\infty}(M,N)\rightarrow\mathbb{R},\quad E_{2}(\varphi)=\frac{1}{2}\int_{M}|\tau(\varphi)|^{2}\ v_g,
$$
and it is characterized by the vanishing of its {\it bitension field}
$$
\tau_2(\varphi)=-\Delta^\varphi \tau(\varphi)-\trace_g R^N(d\varphi,\tau(\varphi))d\varphi,
$$
where
$$
\Delta^\varphi=-\trace_g\left(\nabla^\varphi\nabla^\varphi-\nabla^\varphi_\nabla\right)
$$
is the rough Laplacian of $\varphi^{-1}TN$ and the curvature tensor field is
$$
R^N(X,Y)Z=\nabla^N_X\nabla^N_Y Z-\nabla^N_Y\nabla^N_X Z-\nabla^N_{[X,Y]} Z, \quad \forall X,Y,Z\in C(TM).
$$

We remark that the {\it biharmonic equation} $\tau_2(\varphi)=0$ is a fourth-order non-linear elliptic equation and that any harmonic map is biharmonic. A non-harmonic biharmonic map is called proper biharmonic.

In \cite{J}, G. Y. Jiang defined the stress-energy tensor $S_2$ of the bienergy (also called {\it stress-bienergy tensor}) by
\begin{align*}
  S_2(X,Y) = & \frac{1}{2}|\tau(\varphi)|^2\langle X,Y\rangle +\langle d\varphi,\nabla \tau(\varphi)\rangle \langle X,Y\rangle \\
   & - \langle d\varphi(X),\nabla_Y \tau(\varphi)\rangle - \langle d\varphi(Y),\nabla_X \tau(\varphi)\rangle,
\end{align*}
as it satisfies
$$
\Div S_2=-\langle\tau_2(\varphi),d\varphi\rangle.
$$

The tensor field $S_2$ has a variational meaning, as in the harmonic case. We fix a map $\varphi:M^m\to\left(N^n,h\right)$ and define a new functional
$$
\mathcal{F}_2:\mathcal{G}\to\mathbb{R}, \quad \mathcal{F}_2(g)=E_2(\varphi).
$$
Then we have the following result.
\begin{theorem}[\cite{LMO}]
Let $\varphi:M^m\to\left(N^n,h\right)$ and assume that $M$ is compact. Then
$$
\left.\frac{d}{dt}\right|_{t=0} \mathcal{F}_2\left(g_t\right)=-\frac{1}{2}\int_M \langle \omega,S_2\rangle \ v_g,
$$
so $g$ is a critical point of $\mathcal{F}_2$ if and only if $S_2=0$.
\end{theorem}

We mention that, if $\varphi:\left(M^m,g\right)\to \left(N^n,h\right)$ is an isometric immersion then $\Div S_2$ does not necessarily vanish.

{\it A submanifold of a given Riemannian manifold $\left(N^n,h\right)$} is a pair $\left(M^m,\varphi\right)$, where $M^m$ is a manifold and $\varphi:M\to N$ is an immersion. We always consider on $M$ the induced metric $g=\varphi^\ast h$, thus $\varphi:(M,g)\to (N,h)$ is an isometric immersion; for simplicity we will write $\varphi:M\to N$ without mentioning the metrics. Also, we will write $\varphi:M \to N$, or even $M$, instead of $(M,\varphi)$.

A submanifold $\varphi:M^m\rightarrow N^n$ is called  {\it biharmonic} if the isometric immersion $\varphi$ is a biharmonic map from $\left(M^m,g\right)$ to $\left(N^n,h\right)$.

Even if the notion of biharmonicity may be more appropriate for maps than for submanifolds, as the domain and the codomain metrics are fixed and the variation is made only through the maps, the biharmonic submanifolds proved to be an interesting notion (see, for example, \cite{OH}).

In order to fix the notations, we recall here only the fundamental equations of first order of a submanifold in a Riemannian manifold. These equations define the second fundamental form, the shape operator and the connection in the normal bundle. Let $\varphi:M^m\to N^n$ be an isometric immersion. For each $p\in M$, $T_{\varphi(p)}N$ splits as an orthogonal direct sum
\begin{equation}\label{eq: decomp_immers}
T_{\varphi(p)}N=d\varphi(T_pM)\oplus d\varphi(T_pM)^\perp,
\end{equation}
and $\displaystyle{NM=\bigcup_{p\in M}d\varphi(T_pM)^\perp}$ is referred to as the normal bundle of $\varphi$, or of $M$, in $N$.

Denote by $\nabla$ and $\nabla^N$  the Levi-Civita connections on $M$ and $N$, respectively, and by $\nabla^\varphi$ the induced connection in the pull-back bundle
$\displaystyle{\varphi^{-1}(TN)=\bigcup_{p\in M}T_{\varphi(p)}N}$.
Taking into account the decomposition in \eqref{eq: decomp_immers},
one has
$$
\nabla^\varphi_X d\varphi(Y)=d\varphi(\nabla_X Y)+B(X,Y),\qquad \forall\,X, Y\in C(TM),
$$
where $B\in C(\odot^2 T^\ast M\otimes NM)$ is called the second fundamental form of $M$ in $N$. Here $T^\ast M$ denotes the cotangent bundle of $M$. The mean curvature vector field of $M$ in $N$ is
defined by $H=(\trace B)/m\in C(NM)$, where the $\trace$ is considered with respect to the metric $g$.

Furthermore, if $\eta\in C(NM)$, then
$$
\nabla^\varphi_X \eta=-d\varphi(A_\eta(X))+\nabla^\perp_X\eta, \quad \forall\,X\in C(TM),
$$
where $A_\eta\in C(T^\ast M\otimes TM)$ is called the shape operator of $M$ in $N$ in the direction of $\eta$, and $\nabla^\perp$ is the induced connection in the
normal bundle. Moreover, $\langle B(X,Y),\eta\rangle=\langle A_\eta(X),Y\rangle$, for all $X, Y\in C(TM)$, $\eta\in C(NM)$. In the case of hypersurfaces, we denote $f=\trace A$, where $A=A_\eta$ and $\eta$ is the unit normal vector field, and we have $H=(f/m)\eta$; $f$ is the \textit{($m$ times) mean curvature function}.

A submanifold $M$ of $N$ is called $PMC$ if $H$ is parallel in the normal bundle, and $CMC$ if $|H|$ is constant.

When confusion is unlikely we identify, locally, $M$ with its image through $\varphi$, $X$ with $d\varphi(X)$ and $\nabla^\varphi_X d\varphi(Y)$ with $\nabla^N_X Y$. With these identifications in mind, we write
$$
\nabla^N_X Y=\nabla_X Y+B(X,Y),
$$
and
$$
\nabla^N_X \eta=-A_\eta(X)+\nabla^\perp_X\eta.
$$

If $\Div S_2=0$ for a submanifold $M$ in $N$, then $M$ is called {\it biconservative}. Thus, $M$ is biconservative if and only if the tangent part of its bitension field vanishes.

We have the following characterization theorem of biharmonic submanifolds, obtained by splitting the bitension field in the tangent and normal part.

\begin{theorem}
\label{biharmonic_submanifold}
A submanifold $M^m$ of a Riemannian manifold $N^n$ is biharmonic if and only if
$$
\trace A_{\nabla^\perp_{\cdot} H}(\cdot)+\trace \nabla A_H +\trace \left( R^N(\cdot,H)\cdot\right)^T=0
$$
and
$$
\Delta^\perp H+\trace B\left(\cdot,A_H(\cdot)\right) +\trace \left(R^N(\cdot,H)\cdot\right)^\perp=0,
$$
where $\Delta^\perp$ is the Laplacian in the normal bundle.
\end{theorem}

Various forms of the above result were obtained in \cite{C84,LMO,O02}. From here we deduce some characterization formulas for the biconservativity.

\begin{corollary}
Let $M^m$ be a submanifold of a Riemannian manifold $N^n$. Then $M$ is a biconservative submanifold if and only if:
\begin{enumerate}
    \item $\trace A_{\nabla^\perp_{\cdot} H}(\cdot)+\trace \nabla A_H +\trace \left( R^N(\cdot,H)\cdot\right)^T=0$;
    \item $\frac{m}{2}\grad\left(|H|^2\right)+2\trace A_{\nabla^\perp_{\cdot} H}(\cdot) + 2\trace \left( R^N(\cdot,H)\cdot\right)^T=0$;
    \item $2\trace \nabla A_H-\frac{m}{2}\grad\left(|H|^2\right)=0$.
\end{enumerate}
\end{corollary}

The following properties are immediate.

\begin{proposition}
Let $M^m$ be a submanifold of a Riemannian manifold $N^n$. If $\nabla A_H=0$ then $M$ is biconservative.
\end{proposition}

\begin{proposition}
Let $M^m$ be a submanifold of a Riemannian manifold $N^n$. Assume that $N$ is a space form, i.e., it has constant sectional curvature, and $M$ is $PMC$. Then $M$ is biconservative.
\end{proposition}

\begin{proposition}[\cite{BMO13}]
Let $M^m$ be a submanifold of a Riemannian manifold $N^n$. Assume that $M$ is pseudo-umbilical, i.e., $A_H=|H|^2I$, and $m\neq4$. Then $M$ is $CMC$.
\end{proposition}

If we consider the particular case of hypersurfaces, then Theorem $\ref{biharmonic_submanifold}$ becomes

\begin{theorem}[\cite{BMO13,O10}]
If $M^m$ is a hypersurface in a Riemannian manifold $N^{m+1}$, then $M$ is biharmonic if and only if
$$
2A(\grad f)+f\grad f-2f\left(\ricci^N(\eta)\right)^T=0,
$$
and
$$
\Delta f+f|A|^2-f\ricci^N(\eta,\eta)=0,
$$
where $\eta$ is the unit normal vector field of $M$ in $N$.
\end{theorem}

\begin{corollary}
A hypersurface $M^m$ in a space form $N^{m+1}(c)$ is biconservative if and only if
$$
A(\grad f)=-\frac{f}{2}\grad f.
$$
\end{corollary}

\begin{corollary}
Any $CMC$ hypersurface in $N^{m+1}(c)$ is biconservative.
\end{corollary}

Therefore, the biconservative hypersurfaces may be seen as the next research topic after that of $CMC$ surfaces.

\section{Intrinsic characterization of biconservative surfaces}\label{section2}

We are interested to study biconservative surfaces which are non-$CMC$. We will first look at them from a local, extrinsic point of view and then from a global point of view. While by ``local'' we will mean the biconservative surfaces $\varphi:M^2\to N^3(c)$ with $f>0$ and $\grad f\neq 0$ at any point of $M$, by ``global'' we will mean the \textit{complete} biconservative surfaces $\varphi:M^2\to N^3(c)$ with $f>0$ at any point of $M$ and $\grad f\neq 0$ at any point of an open and dense subset of $M$.

In this section, we consider the local problem, i.e.,  we take $\varphi:M^2\to N^3(c)$ a biconservative surface and assume that $f>0$ and $\grad f\neq 0$ at any point of $M$. Let $X_1=(\grad f)/|\grad f|$ and $X_2$ two vector fields such that $\left\{X_1(p),X_2(p)\right\}$ is a positively oriented orthonormal basis at any point $p\in M$. In particular, we obtain that $M$ is parallelizable. If we denote by $\lambda_1\leq\lambda_2$ the eigenvalues functions of the shape operator $A$, since $A\left(X_1\right)=-(f/2)X_1$ and $\trace A=f$, we get $\lambda_1=-f/2$ and $\lambda_2=3f/2$. Thus the matrix of $A$ with respect to the (global) orthonormal frame field $\left\{X_1,X_2\right\}$ is
\begin{equation*}
A=\left(
\begin{array}{cc}
  -\frac{f}{2} & 0 \\\\
  0 & \frac{3f}{2}
\end{array}
\right).
\end{equation*}
We denote by $K$ the Gaussian curvature and, from the Gauss equation, $K=c+\det A$, we obtain
\begin{equation}\label{eq-f^2}
f^2=\frac{4}{3}(c-K).
\end{equation}
Thus $c-K>0$ on $M$.

From the definitions of $X_1$ and $X_2$, we find that
$$
\grad f=\left(X_1 f\right)X_1 \quad \text{ and } \quad X_2 f=0.
$$
Using the connection $1$-forms, the Codazzi equation and then the extrinsic and intrinsic expression for the Gaussian curvature, we obtain the next result which shows that the mean curvature function of a non-$CMC$ biconservative surface must satisfy a second-order partial differential equation. More precisely, we have the following theorem.

\begin{theorem}[\cite{CMOP}]
Let $\varphi:M^2\to N^3(c)$ a biconservative surface with $f>0$ and $\grad f\neq 0$ at any point of $M$. Then we have
\begin{equation}\label{f-bicons}
f\Delta f+|\grad f|^2+\frac{4}{3}c f^2-f^4=0,
\end{equation}
where $\Delta$ is the Laplace-Beltrami operator on $M$.
\end{theorem}

In fact, we can see that around any point of $M$ there exists $(U;u,v)$ local coordinates such that $f=f(u,v)=f(u)$ and $\eqref{f-bicons}$ is equivalent to
\begin{equation}\label{f-bicons2}
ff^{\prime\prime}-\frac{7}{4}\left(f^\prime\right)^2-\frac{4}{3}cf^2+f^4=0,
\end{equation}
i.e., $f$ must satisfy a \textit{second-order ordinary differential equation}.

Indeed, let $p_0\in M$ be an arbitrary fixed point of $M$ and let $\gamma=\gamma(u)$ be an integral curve of $X_1$ with $\gamma(0)=p_0$. Let $\phi$ the flow of $X_2$ and $(U;u,v)$ local coordinates with $p_0\in U$ such that
$$
X(u,v)=\phi_{\gamma(u)}(v)=\phi(\gamma(u),v).
$$
We have
$$
X_u(u,0)=\gamma^\prime(u)=X_1(\gamma(u))=X_1(u,0)
$$
and
$$
X_v(u,v)=\phi^\prime_{\gamma(u)}(v)=X_2\left(\phi_{\gamma(u)}(v)\right)=X_2(u,v).
$$
If we write the Riemannian metric $g$ on $M$ in local coordinates as
$$
g=g_{11}du^2+2g_{12}du dv+g_{22}dv^2,
$$
we get $g_{22}=\left|X_v\right|^2=\left|X_2\right|^2=1$, and $X_1$ can be expressed with respect to $X_u$ and $X_v$ as
$$
X_1=\frac{1}{\sigma}\left(X_u-g_{12}X_v\right)=\sigma \grad u,
$$
where $\sigma=\sqrt{g_{11}-g_{12}^2}>0$, $\sigma=\sigma(u,v)$.

Let $f\circ X=f(u,v)$. Since $X_2f=0$, we find that
$$
f(u,v)=f(u,0)=f(u), \quad \forall (u,v)\in U.
$$
It can be proved that
$$
\left[X_1,X_2\right]=\frac{3\left(X_1 f\right)}{4f}X_2,
$$
and thus $X_2 X_1 f=X_1 X_2 f-\left[X_1,X_2\right]f=0$.

On the other hand we have
\begin{equation}\label{eqX2X1}
\begin{array}{rl}
  X_2 X_1 f =& X_v\left(\frac{1}{\sigma}f^\prime\right)=X_v\left(\frac{1}{\sigma}\right)f^\prime \\
  = & 0
\end{array}.
\end{equation}
We recall that
$$
\grad f=\left(X_1 f\right)X_1=\left(\frac{1}{\sigma}f^\prime\right)X_1\neq 0
$$
at any point of $U$, and then $f^\prime\neq 0$ at any point of $U$. Therefore, from \eqref{eqX2X1}, $X_v\left(1/\sigma\right)=0$, i.e., $\sigma=\sigma(u)$. Since $g_{11}(u,0)=1$, and $g_{12}(u,0)=0$, we have $\sigma=1$, i.e.,
\begin{equation}\label{eq:X1}
X_1=X_u-g_{12}X_v=\grad u.
\end{equation}
In \cite{CMOP} it was found an equivalent expression for $\eqref{f-bicons}$, i.e.,
$$
\left(X_1X_1f\right)f=\frac{7}{4}\left(X_1 f\right)^2+\frac{4c}{3}f^2-f^4.
$$
Therefore, using \eqref{eq:X1}, relation \eqref{f-bicons} is equivalent to \eqref{f-bicons2}.

\begin{remark}
If $\varphi:M^2\to N^3(c)$ is a non-$CMC$ biharmonic surface, then, there exists an open subset $U$ such that $f>0$, $\grad f\neq 0$ at any point of $U$, and $f$ satisfies the following system
\begin{equation*}
\left\{
  \begin{array}{ll}
    \Delta f=f\left(2c-|A|^2\right) \\\\
    A(\grad f)=-\frac{f}{2}\grad f
  \end{array}
\right. .
\end{equation*}
As we have seen, this system implies
\begin{equation*}
\left\{
  \begin{array}{ll}
    \Delta f=f\left(2c-|A|^2\right) \\\\
    f\Delta f+|\grad f|^2+\frac{4}{3}c f^2-f^4=0
  \end{array}
\right. .
\end{equation*}
which, in fact, is a ODE system. We get
\begin{equation}\label{sis}
\left\{
  \begin{array}{ll}
    ff^{\prime\prime}-\frac{3}{4}\left(f^\prime\right)^2+2cf^2-\frac{5}{2}f^4=0\\\\
    ff^{\prime\prime}-\frac{7}{4}\left(f^\prime\right)^2-\frac{4}{3}cf^2+f^4=0
  \end{array}
\right. .
\end{equation}
As an immediate consequence we obtain
$$
\left(f^{\prime}\right)^2+\frac{10}{3}cf^2-\frac{7}{2}f^4=0,
$$
and combining it with the first integral
$$
\left(f^\prime\right)^2=2f^4-8cf^2+\alpha f^{3/2}
$$
of the first equation from \eqref{sis}, where $\alpha\in\mathbb{R}$ is a constant, we obtain
$$
\frac{3}{2}f^{5/2}+\frac{14}{3}cf^{1/2}-\alpha=0.
$$
If we denote $\tilde{f}=f^{1/2}$, we get ${3\tilde{f}^{5}}/{2}+{14 c\tilde{f}}/{3}-\alpha=0$. Thus, $\tilde{f}$ satisfies a polynomial equation with constant coefficients, so $\tilde{f}$ has to be a constant and then, $f$ is a constant, i.e., $\grad f=0$ on $U$ (in fact, $f$ has to be zero). Therefore, we have a contradiction (see \cite{C91,CI91} for $c=0$ and \cite{CMO02,CMO01}, for $c=\pm 1$).
\end{remark}

We can also note that relation $\eqref{f-bicons}$, which is an extrinsic relation, together with \eqref{eq-f^2}, allows us to find an \textit{intrinsic relation} that $(M,g)$ must satisfy. More precisely, the Gaussian curvature of $M$ has to satisfy
\begin{equation}\label{bicons-K}
(c-K)\Delta K-|\grad K|^2-\frac{8}{3}K(c-K)^2=0,
\end{equation}
and the conditions $c-K>0$ and $\grad K\neq0$.

Formula \eqref{bicons-K} is very similar to the Ricci condition. Further, we will briefly recall the Ricci problem. Given an abstract surface $\left(M^2,g\right)$, we want to find the conditions that have to be satisfied by $M$ such that, locally, it admits a minimal embedding in $N^3(c)$. It was proved (see \cite{MM,R}) that if $\left(M^2,g\right)$ is an abstract surface such that $c-K>0$ at any point of $M$, where $c\in \mathbb{R}$ is a constant, then, locally, it admits a minimal embedding in $N^3(c)$ if and only if
\begin{equation}\label{Ricci}
(c-K)\Delta K-|\grad K|^2-4K(c-K)^2=0.
\end{equation}
Condition \eqref{Ricci} is called the \textit{Ricci condition with respect to $c$}, or simply the \textit{Ricci condition}. If $\eqref{Ricci}$ holds, then, locally, $M$ admits a one-parameter family of minimal embeddings in $N^3(c)$.

We can see that relations $(\ref{bicons-K})$ and $\eqref{Ricci}$ are very similar and, in \cite{FNO}, the authors studied the link between them. Thus, for $c=0$, it was proved that if we consider a surface $\left(M^2,g\right)$ which satisfies $\eqref{bicons-K}$ and $K<0$, then there exists a very simple conformal transformation of the metric $g$ such that $\left(M^2,\sqrt{-K}g\right)$ satisfies $\eqref{Ricci}$. A similar result was also proved for $c\neq 0$, but in this case, the conformal factor has a complicated expression (and it is not enough to impose that $\left(M^2,g\right)$ satisfy \eqref{bicons-K}, but we need the stronger hypothesis of it to admit a non-$CMC$ biconservative immersion in $N^3(c)$).

Unfortunately, condition $\eqref{bicons-K}$ does not imply, locally, the existence of a biconservative immersion in $N^3(c)$, as in the minimal case. We need a stronger condition. It was obtained the following local, intrinsic characterization theorem.

\begin{theorem}[\cite{FNO}]\label{thm:char}
Let $\left(M^2,g\right)$ be an abstract surface and $c\in\mathbb{R}$ a constant. Then, locally, $M$ can be isometrically embedded in a space form $N^3(c)$ as a biconservative surface with positive mean curvature having the gradient different from zero at any point if and only if the Gaussian curvature $K$ satisfies $c-K(p)>0$, $(\grad K)(p)\neq 0$, for any point $p\in M$, and its level curves are circles in $M$ with constant curvature
$$
\kappa=\frac{3|\grad K|}{8(c-K)}.
$$
\end{theorem}

\begin{remark}
If the surface $M$ in Theorem \ref{thm:char} is simply connected, then the theorem holds globally, but, in this case, instead of a local isometric embedding we have a global isometric immersion.
\end{remark}

We remark that unlike in the minimal immersions case, if $M$ satisfies the hypotheses from Theorem \ref{thm:char}, then there exists a \textit{unique} biconservative immersion in $N^3(c)$ (up to an isometry of $N^3(c)$), and not a one-parameter family.

The characterization theorem can be equivalently rewritten as below.

\begin{theorem}\label{thm:char-local}
Let $\left(M^2,g\right)$ be an abstract surface with Gaussian curvature $K$ satisfying $c-K(p)>0$ and $(\grad K)(p)\neq 0$ at any point $p\in M$, where $c\in \mathbb{R}$ is a constant. Let $X_1=(\grad K)/|\grad K|$ and $X_2\in C(TM)$ be two vector fields on $M$ such that $\left\{X_1(p),X_2(p)\right\}$ is a positively oriented basis at any point of $p\in M$. Then, the following conditions are equivalent:
\begin{itemize}
  \item [(a)] the level curves of $K$ are circles in $M$ with constant curvature
  $$
  \kappa=\frac{3|\grad K|}{8(c-K)}=\frac{3X_1K}{8(c-K)};
  $$
  \item [(b)]
  $$
  X_2\left(X_1K\right)=0\quad \text{and} \quad \nabla_{X_2}X_2=\frac{-3X_1K}{8(c-K)}X_1;
  $$
  \item [(c)] locally, the metric $g$ can be written as $g=(c-K)^{-3/4}\left(du^2+dv^2\right)$, where $(u,v)$ are local coordinates positively oriented, $K=K(u)$, and $K^\prime>0$;
  \item [(d)] locally, the metric $g$ can be written as $g=e^{2\varphi}\left(du^2+dv^2\right)$, where $(u,v)$ are local coordinates positively oriented, and $\varphi=\varphi(u)$ satisfies the equation
      \begin{equation}\label{ec-d}
      \varphi^{\prime\prime}=e^{-2\varphi/3}-ce^{2\varphi}
      \end{equation}
      and the condition $\varphi^\prime>0$; moreover, the solutions of the above equation, $u=u(\varphi)$, are
      $$
      u=\bigintsss_{\varphi_0}^{\varphi}\frac{d\tau}{\sqrt{-3e^{-2\tau/3}-ce^{2\tau}+a}}+u_0,
      $$
      where $\varphi$ is in some open interval $I$ and $a,u_0\in \mathbb{R}$ are constants;
  \item [(e)] locally, the metric $g$ can be written as $g=e^{2\varphi}\left(du^2+dv^2\right)$, where $(u,v)$ are local coordinates positively oriented, and $\varphi=\varphi(u)$ satisfies the equation
      \begin{equation}\label{ec-e}
      3\varphi^{\prime\prime\prime}+2\varphi^\prime\varphi^{\prime\prime}+8ce^{2\varphi}\varphi^\prime=0
      \end{equation}
      and the conditions $\varphi^\prime>0$ and $c+e^{-2\varphi}\varphi^{\prime\prime}>0$; moreover, the solutions of the above equation, $u=u(\varphi)$, are
      $$
      u=\bigintsss_{\varphi_0}^{\varphi}\frac{d\tau}{\sqrt{-3be^{-2\tau/3}-ce^{2\tau}+a}}+u_0,
      $$
      where $\varphi$ is in some open interval $I$ and $a,b,u_0\in \mathbb{R}$ are constants, $b>0$.
\end{itemize}
\end{theorem}

The proof follows by direct computations and by using Remark 4.3 in \cite{FNO} and Proposition 3.4 in \cite{N}.
\begin{remark}
From the above theorem we have the following remarks.
\begin{itemize}
  \item [(i)] If condition $(a)$ is satisfied, i.e., the integral curves of $X_2$ are circles in $M$ with a precise constant curvature, then the integral curves of $X_1$ are geodesics of $M$.

  \item [(ii)] If condition $(c)$ is satisfied, then $K$ has to be a solution of the equation
      $$
      3K^{\prime\prime}(c-K)+3\left(K^\prime\right)^2+8K(c-K)^{5/4}=0.
      $$

  \item [(iii)] If condition $(c)$ is satisfied and $c>0$, then  $\left(M^2, (c-K)^{3/4}g\right)$ is a flat surface and, trivially, a Ricci surface with respect to $c$.

  \item [(iv)] Let $\varphi=\varphi(u)$ be a solution of equation \eqref{ec-e}. We consider the change of coordinates
      $$
      (u,v)=\left(\alpha\tilde{u}+\beta,\alpha\tilde{v}+\beta\right),
      $$
      where $\alpha\in \mathbb{R}$ is a positive constant and $\beta\in\mathbb{R}$, and define
      $$
      \phi=\varphi\left(\alpha\tilde{u}+\beta\right)+\log \alpha.
      $$
      Then $g=e^{2\phi}\left(d\tilde{u}^2+d\tilde{v}^2\right)$ and $\phi$ also satisfies equation \eqref{ec-e}. If $\varphi=\varphi(u)$ satisfies the first integral
      $$
      \varphi^{\prime\prime}=be^{-2\varphi/3}-ce^{2\varphi},
      $$
      where $b>0$, then, for $\alpha=b^{-3/8}$, $\phi=\phi\left(\tilde{u}\right)$ satisfies
      $$
      \phi^{\prime\prime}=e^{-2\phi/3}-ce^{2\phi}.
      $$
      From here, as the classification is done up to isometries, we note that the parameter $b$ in the solution of \eqref{ec-e} is not essential and only the parameter $a$ counts. Thus we have a one-parameter family of solutions.

  \item [(v)] If $\varphi$ is a solution of \eqref{ec-e}, for some $c$, then $\varphi+\alpha$, where $\alpha$ is a real constant, is a solution of \eqref{ec-e} for $ce^{2\alpha}$.

  \item [(vi)] If $c=0$, we note that if $\varphi$ is a solution of \eqref{ec-e}, then also $\varphi+\cst$ is a solution of the same equation, i.e, condition $(a)$ from Theorem $\ref{thm:char-local}$ is invariant under the homothetic tranformations of the metric $g$. Then, we see that equation \eqref{ec-e} is invariant under the affine change of parameter $u=\alpha\tilde{u}+\beta$, where $\alpha>0$. Therefore, we must solve equation \eqref{ec-e} up to this change of parameter and an additive constant of the solution $\varphi$. The additive constant will be the parameter that counts.
\end{itemize}
\end{remark}

In the $c=0$ case, the solutions of equation \eqref{ec-e}, are explicitly determined in the next proposition.

\begin{proposition}[\cite{N}]\label{prop-c=0}
The solutions of the equation
$$
3\varphi^{\prime\prime\prime}+2\varphi^\prime\varphi^{\prime\prime}=0
$$
which satisfy the conditions $\varphi^\prime>0$ and $\varphi^{\prime\prime}>0$, up to affine transformations of the parameter with $\alpha>0$, are given by
$$
\varphi(u)=3\log(\cosh u)+\cst, \qquad u>0.
$$
\end{proposition}

We note that, when $c=0$, we have a one-parameter family of solutions of equation \eqref{ec-e}, i.e., $g_{C_0}=C_0(\cosh u)^6\left(du^2+dv^2\right)$, $C_0$ being a positive constant.

If $c\neq0$, then we can not determine explicitly $\varphi=\varphi(u)$. Another way to see that in the $c\neq 0$ case we have only a one-parameter family of solutions of equation \eqref{ec-e} is to rewrite the metric $g$ in certain non-isothermal coordinates.

Further, we will consider only the $c=1$ case and we have the next result.

\begin{proposition}[\cite{N}]
Let $\left(M^2,g\right)$ be an abstract surface with $g=e^{2\varphi(u)}(du^2+dv^2)$, where $u=u(\varphi)$ satisfies
$$
u=\bigintsss_{\varphi_0}^{\varphi}\frac{d\tau}{\sqrt{-3be^{-2\tau/3}-e^{2\tau}+a}}+u_0,
$$
where $\varphi$ is in some open interval $I$, $a,b\in \mathbb{R}$ are positive constants, and $u_0\in\mathbb{R}$ is a constant. Then $\left(M^2,g\right)$ is isometric to
$$
\left(D_{C_1},g_{C_1}=\frac{3}{\xi^2\left(-\xi^{8/3}+3C_1\xi^2-3\right)}d\xi^2 +\frac{1}{\xi^2}d\theta^2\right),
$$
where $D_{C_1}=\left(\xi_{01},\xi_{02}\right)\times\mathbb{R}$, $C_1\in\left(4/\left(3^{3/2}\right),\infty\right)$ is a positive constant, and $\xi_{01}$ and $\xi_{02}$ are the positive vanishing points of $-\xi^{8/3}+3C_1\xi^2-3$, with $0<\xi_{01}<\xi_{02}$.
\end{proposition}

\begin{remark}
Let us consider
$$
\left(D_{C_1},g_{C_1}=\frac{3}{\xi^2\left(-\xi^{8/3}+3{C_1}\xi^2-3\right)}d\xi^2+ \frac{1}{\xi^2}d\theta^2\right)
$$
and
$$
\left(D_{C_1^\prime},g_{C_1^\prime}=\frac{3}{\tilde{\xi}^2\left(-\tilde{\xi}^{8/3}+3C_1^\prime \tilde{\xi}^2-3\right)}d\tilde{\xi}^2+\frac{1}{\tilde{\xi}^2}d\tilde{\theta}^2\right).
$$
The surfaces $\left(D_{C_1},g_{C_1}\right)$ and $\left(D_{C_1^\prime},g_{C_1^\prime}\right)$ are isometric if and only if $C_1=C_1^\prime$ and the isometry is $\Theta(\xi,\theta)=\left(\xi,\pm\theta+\cst\right)$. Therefore, we have a one-parameter family of surfaces.
\end{remark}

\begin{remark}
We note that the expression of the Gaussian curvature of $\left(D_{C_1},g_{C_1}\right)$ does not depend on $C_1$. More precisely,
$$
K_{C_1}(\xi,\theta)=-\frac{1}{9}\xi^{8/3}+1.
$$
But, if we change further the coordinates $(\xi,\theta)=\left(\xi_{01}+\tilde{\xi}\left(\xi_{02}-\xi_{01}\right),\tilde{\theta}\right)$, then we ``fix'' the domain, i.e., $\left(D_{C_1},g_{C_1}\right)$ is isometric to $\left((0,1),\tilde{g}_{C_1}\right)$ and $C_1$ appears in the expression of $K_{C_1}\left(\tilde{\xi},\tilde{\theta}\right)$.
\end{remark}

\section{Complete biconservative surfaces in $\mathbb{R}^3$}

In this section we consider the global problem and construct complete biconservative surfaces in $\mathbb{R}^3$ with $f>0$ everywhere and $\grad f\neq 0$ at any point of an open dense subset. Or, from intrinsic point of view, we construct a complete abstract surface $\left(M^2,g\right)$ with $K<0$ everywhere and $\grad K\neq 0$ at any point of an open dense subset of $M$, that admits a biconservative immersion in $\mathbb{R}^3$, defined on the whole $M$, with $f>0$ on $M$ and $|\grad f|>0$ on the open dense subset.

First, we recall a local extrinsic result which provides a characterization of biconservative surfaces in $\mathbb{R}^3$.

\begin{theorem}[\cite{HV95}]
Let $M^2$ be a surface in $\mathbb{R}^3$ with $f(p)>0$ and $(\grad f)(p)\neq0$ for any $p\in M$. Then, $M$ is biconservative if and only if, locally, it is a surface of revolution, and the curvature $\kappa=\kappa(u)$ of the profile curve $\sigma=\sigma(u)$, $\left|\sigma^\prime(u)\right|=1$, is a positive solution of the following ODE
$$
\kappa^{\prime\prime}\kappa=\frac{7}{4}\left(\kappa^\prime\right)^2-4\kappa^4.
$$
\end{theorem}

In \cite{CMOP} there was found the local explicit parametric equation of a biconservative surface in $\mathbb{R}^3$.

\begin{theorem}[\cite{CMOP}]
Let $M^2$ be a biconservative surface in $\mathbb{R}^3$ with $f(p)>0$ and $(\grad f)(p)\neq0$ for any $p\in M$. Then, locally, the surface can be parametrized by
$$
X_{\tilde{C}_0}(\rho,v)=\left(\rho\cos v,\rho \sin v, u_{\tilde{C}_0}(\rho)\right),
$$
where
$$
u_{\tilde{C}_0}(\rho)=\frac{3}{2\tilde{C}_0}\left(\rho^{1/3}\sqrt{\tilde{C}_0\rho^{2/3}-1}+\frac{1}{\sqrt{\tilde{C}_0}}\log\left(\sqrt{\tilde{C}_0}\rho^{1/3}+\sqrt{\tilde{C}_0\rho^{2/3}-1}\right)\right)
$$
with $\tilde{C}_0$ a positive constant and $\rho\in\left({\tilde{C}_0}^{-3/2},\infty\right)$.
\end{theorem}

We denote by $S_{\tilde{C}_0}$ the image $X_{\tilde{C}_0}\left(\left({\tilde{C}_0}^{-3/2},\infty\right)\times\mathbb{R}\right)$. We note that any two such surfaces are not locally isometric, so we have a one-parameter family of biconservative surfaces in $\mathbb{R}^3$. These surfaces are not complete.

\begin{remark}
If $\varphi:M^2\to\mathbb{R}^3$ is a biconservative surface with $f>0$ and $\grad f\neq 0$ at any point, then there exists a unique $\tilde{C}_0$ such that $\varphi(M)\subset S_{\tilde{C}_0}$. Indeed, any point admits an open neighborhood which is an open subset of some $S_{\tilde{C}_0}$. Let us consider $p_0\in M$. Then, there exists a unique $\tilde{C}_0$ such that $\varphi(U)\subset S_{\tilde{C}_0}$, where $U$ is an open neighborhood of $p_0$. If $A$ denotes the set of all points of $M$ such that they admit open neighborhoods which are open subsets of that $S_{\tilde{C}_0}$, then the set $A$ is non-empty, open and closed in $M$. Thus, as $M$ is connected, it follows that $A=M$.
\end{remark}

The ``boundary'' of $S_{\tilde{C}_0}$, i.e., $\overline{S}_{\tilde{C}_0}\setminus S_{\tilde{C}_0}$, is the circle $\left({\tilde{C}_0}^{-3/2}\cos v,{\tilde{C}_0}^{-3/2}\sin v,0\right)$, which lies in the $Oxy$ plane. At a boundary point, the tangent plane to the closure $\overline{S}_{\tilde{C}_0}$ of $S_{\tilde{C}_0}$ is parallel to $Oz$. Moreover, along the boundary, the mean curvature function is constant $f_{\tilde{C}_0}=\left(2{\tilde{C}_0}^{3/2}\right)/3$ and $\grad f_{\tilde{C}_0}=0$.

Thus, in order to obtain a complete biconservative surface in $\mathbb{R}^3$, we can expect to ``glue'' along the boundary two biconservative surfaces of type $S_{\tilde{C}_0}$ corresponding to the same $\tilde{C}_0$ (the two constants have to be the same) and symmetric to each other, at the level of $C^\infty$ smoothness.

In fact, it was proved that we can glue two biconservative surfaces $S_{\tilde{C}_0}$ and $S_{{\tilde{C}_0}^\prime}$, at the level of $C^\infty$ smoothness, only along the boundary and, in this case, $\tilde{C}_0={\tilde{C}_0}^\prime$.

\begin{proposition}[\cite{MOR16,N}]
If we consider the symmetry of $\Graf u_C$, with respect to the $O\rho(=Ox)$ axis, we get a smooth, complete, biconservative surface $\tilde{S}_{\tilde{C}_0}$ in $\mathbb{R}^3$. Moreover, its mean curvature function $\tilde{f}_{\tilde{C}_0}$ is positive and $\grad \tilde{f}_{\tilde{C}_0}$ is different from zero at any point of an open dense subset of $\tilde{S}_{\tilde{C}_0}$.
\end{proposition}

\begin{remark}
The profile curve $\sigma_{\tilde{C}_0}=\left(\rho,0,u_{\tilde{C}_0}(\rho)\right)\equiv \left(\rho,u_{\tilde{C}_0}(\rho)\right)$ can be re\-para\-metrized as
\begin{equation}\label{eq:sigma_tc}
\begin{array}{rl}
\sigma_{\tilde{C}_0}(\theta)= & \left(\sigma_{\tilde{C}_0}^1(\theta), \sigma_{\tilde{C}_0}^2(\theta)\right) \\\\
  = & {\tilde{C}_0}^{-3/2}\left((\theta+1)^{3/2},\frac{3}{2}\left(\sqrt{\theta^2+\theta} + \log \left(\sqrt{\theta}+\sqrt{\theta+1}\right)\right)\right), \qquad \theta>0,
\end{array}
\end{equation}
and now $X_{\tilde{C}_0}=X_{\tilde{C}_0}(\theta,v)$.
\end{remark}

\begin{proposition}
The homothety of $\mathbb{R}^3$, $(x,y,z)\to \tilde{C}_0(x,y,z)$, renders $\tilde{S}_1$ onto $\tilde{S}_{{\tilde{C}_0}^{-2/3}}$.
\end{proposition}

In \cite{N}, there were also found the complete biconservative surfaces in $\mathbb{R}^3$ with $f>0$ at any point and $\grad f\neq 0$ at any point of an open dense subset, but there, the idea was to use the intrinsic characterization of the biconservative surfaces. More precisely, we have the next global result.

\begin{theorem}[\cite{N}]\label{main_th1}
Let $\left(\mathbb{R}^2,g_{C_0}=C_0 \left(\cosh u\right)^6\left(du^2+dv^2\right)\right)$ be a surface, where $C_0\in\mathbb{R}$ is a positive constant. Then we have:
\begin{itemize}
\item[(a)] the metric on $\mathbb{R}^2$ is complete;
\item[(b)] the Gaussian curvature is given by
           $$
            K_{C_0}(u,v)=K_{C_0}(u)=-\frac{3}{C_0\left(\cosh u\right)^8}<0,\quad  K^\prime_{C_0}(u)=\frac{24 \sinh u}{C_0\left(\cosh u\right)^9},
           $$
           and therefore $\grad K_{C_0}\neq 0$ at any point of $\mathbb{R}^2\setminus Ov$;

\item[(c)] the immersion $\varphi_{C_0}:\left(\mathbb{R}^2,g_{C_0}\right)\to \mathbb{R}^3$ given by
    $$
    \varphi_{C_0}(u,v)=\left(\sigma_{C_0}^1(u)\cos (3v), \sigma_{C_0}^1(u)\sin (3v), \sigma_{C_0}^2(u)\right)
    $$
    is biconservative in $\mathbb{R}^3$, where
    $$
    \sigma_{C_0}^1(u)=\frac{\sqrt{C_0}}{3}\left(\cosh u\right)^3, \quad
    \sigma_{C_0}^2(u)=\frac{\sqrt{C_0}}{2}\left(\frac{1}{2}\sinh (2u)+u\right), \qquad u \in \mathbb{R}.
    $$
\end{itemize}
\end{theorem}

\begin{sproof}
The first two items follow by standard arguments. For the last part, we note that choosing $\tilde{C}_0=(9/C_0)^{1/3}$ in  \eqref{eq:sigma_tc} and using the change of coordinates $(\theta,v)=\left((\sinh u)^2,3v\right)$, where $u>0$, the metric induced by $X_{(9/C_0)^{1/3}}$ coincides with $g_{C_0}$. Then, we define $\varphi_{C_0}$ as: for $u>0$, $\varphi_{C_0}(u,v)$ is obtained by rotating the profile curve
$$
\sigma^+_{\left(\frac{9}{C_0}\right)^{1/3}}(u)=\sigma_{\left(\frac{9}{C_0}\right)^{1/3}}(u)=\left(\sigma_{\left(\frac{9}{C_0}\right)^{1/3}}^1(u), \sigma_{\left(\frac{9}{C_0}\right)^{1/3}}^2(u)\right),
$$
and for $u<0$, $\varphi_{C_0}(u,v)$ is obtained by rotating the profile curve
$$
\sigma^-_{\left(\frac{9}{C_0}\right)^{1/3}}(u)=\left(\sigma_{\left(\frac{9}{C_0}\right)^{1/3}}^1(-u), -\sigma_{\left(\frac{9}{C_0}\right)^{1/3}}^2(-u)\right).
$$
\end{sproof}

By simple transformations of the metric, $\left(\mathbb{R}^2,g_{C_0}\right)$ becomes a Ricci surface or a surface with constant Gaussian curvature.

\begin{theorem}
Consider the surface $\left(\mathbb{R}^2,g_{C_0}\right)$. Then $\left(\mathbb{R}^2,\sqrt{-K_{C_0}}g_{C_0}\right)$ is complete, satisfies the Ricci condition and can be minimally immersed in $\mathbb{R}^3$ as a helicoid or a catenoid.
\end{theorem}

\begin{proposition}
Consider the surface $\left(\mathbb{R}^2,g_{C_0}\right)$. Then $\left(\mathbb{R}^2,-K_{C_0}g_{C_0}\right)$ has constant Gaussian curvature $1/3$ and it is not complete. Moreover, $\left(\mathbb{R}^2,-K_{C_0}g_{C_0}\right)$ is the universal cover of the surface of revolution in $\mathbb{R}^3$ given by
$$
Z(u,v)=\left(\alpha(u)\cosh\left(\frac{\sqrt{3}}{a}v\right),\alpha(u)\sinh\left(\frac{\sqrt{3}}{a}v\right),\beta(u)\right), \qquad (u,v)\in\mathbb{R}^2,
$$
where $a\in (0,\sqrt{3}]$ and
$$
\alpha(u)=\frac{a}{\cosh u}, \ \beta(u)=\bigintsss_{0}^{u}{\frac{\sqrt{\left(3-a^2\right)\cosh^2\tau+a^2}}{\cosh^2\tau}\ d\tau}.
$$
\end{proposition}

\begin{remark}
When $a=\sqrt{3}$, the immersion $Z$ has only umbilical points and the image $Z\left(\mathbb{R}^2\right)$ is the round sphere of radius $\sqrt{3}$, without the North and the South poles. Moreover, if $a\in (0,\sqrt{3})$, then $Z$ has no umbilical points.
\end{remark}

\vspace{0.5cm}

Concerning the biharmonic surfaces in $\mathbb{R}^3$ we have the following non-existence result.

\begin{theorem}[\cite{C91,CI91}]
There exists no proper biharmonic surface in $\mathbb{R}^3$.
\end{theorem}

\section{Complete biconservative surfaces in $\mathbb{S}^3$}

As in the previous section, we consider the global problem for biconservative surfaces in $\mathbb{S}^3$, i.e., our aim is to construct complete biconservative surfaces in $\mathbb{S}^3$ with $f>0$ everywhere and $\grad f\neq 0$ at any point of an open and dense subset.

We start with the following local, extrinsic result.

\begin{theorem}[\cite{CMOP}]
Let $M^2$ be a biconservative surface in $\mathbb{S}^3$ with $f(p)>0$ and $(\grad f)(p)\neq 0$ at any point $p\in M$. Then, locally, the surface, viewed in $\mathbb{R}^4$, can be parametrized by
\begin{equation*}
Y_{\tilde{C}_1}(u,v)=\sigma(u)+\frac{4\kappa(u)^{-3/4}}{3\sqrt{\tilde{C}_1}}\left( \overline{f}_1 (\cos v -1)+\overline{f}_2 \sin v\right),
  \end{equation*}
where $\tilde{C}_1\in\left(64/\left(3^{5/4}\right),\infty\right)$ is a positive constant; $\overline{f}_1, \overline{f}_2\in\mathbb{R}^4$ are two constant orthonormal vectors; $\sigma(u)$ is a curve parametrized by arclength that satisfies
\begin{equation*}
  \langle \sigma(u),\overline{f}_1\rangle = \frac{4\kappa(u)^{-3/4}}{3\sqrt{\tilde{C}_1}}, \qquad \langle \sigma(u),\overline{f}_2\rangle=0,
\end{equation*}
and, as a curve in $\mathbb{S}^2$, its curvature $\kappa=\kappa(u)$ is a positive non constant solution of the following ODE
$$
\kappa^{\prime\prime}\kappa=\frac{7}{4}\left(\kappa^\prime\right)^2+\frac{4}{3}\kappa^2-4\kappa^4
$$
such that
$$
\left(\kappa^\prime\right)^2=-\frac{16}{9}\kappa^2-16\kappa^4+\tilde{C}_1\kappa^{7/2}.
$$
\end{theorem}

\begin{remark}
The constant $\tilde{C}_1$ determines uniquely the curvature $\kappa$, up to a translation of $u$, and then $\kappa$, $\overline{f}_1$ and $\overline{f}_2$ determines uniquely the curve $\sigma$.
\end{remark}

We consider $\overline{f}_1=\overline{e}_3$ and $\overline{f}_2=\overline{e}_4$ and change the coordinates $(u,v)$ in $(\kappa,v)$. Then, we get
\begin{equation}\label{eq:Y_{C_1}}
\begin{array}{rl}
  Y_{\tilde{C}_1}(\kappa,v)= & \Bigg(\sqrt{1-\left(\frac{4}{3\sqrt{\tilde{C}_1}}\kappa^{-3/4}\right)^2}\cos \mu(\kappa), \sqrt{1-\left(\frac{4}{3\sqrt{\tilde{C}_1}}\kappa^{-3/4}\right)^2}\sin \mu(\kappa), \\
   & \frac{4}{3\sqrt{\tilde{C}_1}}\kappa^{-3/4}\cos v, \frac{4}{3\sqrt{\tilde{C}_1}}\kappa^{-3/4}\sin v\Bigg),
\end{array}
\end{equation}
where $(\kappa,v)\in \left(\kappa_{01},\kappa_{02}\right)\times\mathbb{R}$, $\kappa_{01}$ and $\kappa_{02}$ are positive solutions of
$$
-\frac{16}{9}\kappa^2-16\kappa^4+\tilde{C}_1\kappa^{7/2}=0
$$
and
$$
\mu(\kappa)=\pm 108\bigintsss_{\kappa_0}^{\kappa}\frac{\sqrt{\tilde{C}_1}\tau^{3/4}}{\left(-16+9\tilde{C}_1\tau^{3/2}\right) \sqrt{ 9\tilde{C}_1\tau^{3/2}-16\left(1+9\tau^2\right)}}\ d\tau +c_0,
$$
with $c_0\in\mathbb{R}$ a constant and $\kappa_0\in\left(\kappa_{01},\kappa_{02}\right)$. We note that an alternative expression for $Y_{\tilde{C}_1}$ was given in \cite{Fu}.

\begin{remark}
The limits $\lim_{\kappa\searrow \kappa_{01}} \mu(\kappa)=\mu\left(\kappa_{01}\right)$ and  $\lim_{\kappa\nearrow \kappa_{02}}\mu(\kappa)=\mu\left(\kappa_{02}\right)$ are finite.
\end{remark}

\begin{remark}
For simplicity, we choose $\kappa_0=(3\tilde{C}_1/64)^2$.
\end{remark}

If we denote $S_{\tilde{C}_1}$ the image of $Y_{\tilde{C}_1}$, then we note that the boundary of $S_{\tilde{C}_1}$ is made up from two circles and along the boundary, the mean curvature function is constant (two different constants) and its gradient vanishes. More precisely, the boundary of $S_{\tilde{C}_1}$ is given by the curves
\begin{equation*}
\begin{array}{c}
\Bigg(\sqrt{1-\left(\frac{4}{3\sqrt{\tilde{C}_1}}\kappa_{01}^{-3/4}\right)^2}\cos\mu\left(\kappa_{01}\right), \sqrt{1-\left(\frac{4}{3\sqrt{\tilde{C}_1}}\kappa_{01}^{-3/4}\right)^2}\sin\mu\left(\kappa_{01}\right), \\
\frac{4}{3\sqrt{\tilde{C}_1}}\kappa_{01}^{-3/4}\cos v, \frac{4}{3\sqrt{\tilde{C}_1}}\kappa_{01}^{-3/4}\sin v \Bigg)
\end{array}
\end{equation*}
and
\begin{equation*}
\begin{array}{c}
\Bigg(\sqrt{1-\left(\frac{4}{3\sqrt{\tilde{C}_1}}\kappa_{02}^{-3/4}\right)^2}\cos\mu\left(\kappa_{02}\right), \sqrt{1-\left(\frac{4}{3\sqrt{\tilde{C}_1}}\kappa_{02}^{-3/4}\right)^2}\sin\mu\left(\kappa_{02}\right), \\
\frac{4}{3\sqrt{\tilde{C}_1}}\kappa_{02}^{-3/4}\cos v, \frac{4}{3\sqrt{\tilde{C}_1}}\kappa_{02}^{-3/4}\sin v \Bigg).
\end{array}
\end{equation*}
These curves are circles in affine planes in $\mathbb{R}^4$ parallel to the $Ox^3x^4$ plane and
their radii are $\left(4\kappa_{01}^{-3/4}\right)/\left(3\sqrt{\tilde{C}_1}\right)$ and $\left(4\kappa_{02}^{-3/4}\right)/\left(3\sqrt{\tilde{C}_1}\right)$, respectively.

At a boundary point, using the coordinates $(\mu, v)$, we get that the tangent plane to the closure of $S_{\tilde{C}_1}$ is spanned by a vector which is tangent to the corresponding circle and by
\begin{equation*}
\begin{array}{c}
\left(-\sqrt{1-\left(\frac{4}{3\sqrt{\tilde{C}_1}}\kappa_{0i}^{-3/4}\right)^2}\sin\mu\left(\kappa_{0i}\right), \sqrt{1-\left(\frac{4}{3\sqrt{\tilde{C}_1}}\kappa_{0i}^{-3/4}\right)^2}\cos\mu\left(\kappa_{0i}\right),0,0 \right),
\end{array}
\end{equation*}
where $i=1$ or $i=2$.

Thus, in order to construct a complete biconservative surface in $\mathbb{S}^3$, we can expect to glue along the boundary two biconservative surfaces of type $S_{\tilde{C}_1}$, corresponding to the same $\tilde{C}_1$. In fact, if we want to glue two surfaces corresponding to $\tilde{C}_1$ and $\tilde{C}_1^\prime$ along the boundary, then these constants have to coincide and there is no ambiguity concerning along which circle of the boundary we should glue the two pieces. But this process is not as clear as in $\mathbb{R}^3$ since we should repeat it infinitely many times.

Further, as in the $\mathbb{R}^3$ case, we change the point of view and use the intrinsic characterization of the biconservative surfaces in $\mathbb{S}^3$.

The surface $\left(D_{C_1},g_{C_1}\right)$ defined in Section \ref{section2} is not complete but it has the following properties.

\begin{theorem}[\cite{N}]
Consider $\left(D_{C_1}, g_{C_1}\right)$. Then, we have
\begin{itemize}
  \item[(a)] $K_{C_1}(\xi,\theta)=K(\xi,\theta)$,
  $$
  1-K(\xi,\theta)=\frac{1}{9}\xi^{8/3}>0,\quad K^\prime(\xi)=-\frac{8}{27}\xi^{5/3}
  $$
  and $\grad K\neq 0$ at any point of $D_{C_1}$;
  \item[(b)] the immersion $\phi_{C_1}:\left(D_{C_1}, g_{C_1}\right)\to \mathbb{S}^3$ given by
  $$
  \phi_{C_1}(\xi,\theta)=\left(\sqrt{1-\frac{1}{C_1\xi^2}}\cos \zeta(\xi),\sqrt{1-\frac{1}{C_1\xi^2}}\sin \zeta(\xi),\frac{\cos(\sqrt{C_1}\theta)}{\sqrt{C_1}\xi}, \frac{\sin(\sqrt{C_1}\theta)}{\sqrt{C_1}\xi}\right),
  $$
  is biconservative in $\mathbb{S}^3$, where
  $$
  \zeta(\xi)=\pm\bigintsss_{\xi_{00}}^\xi{\frac{\sqrt{C_1}\tau^{4/3}}{(-1+C_1\tau^2) \sqrt{-\tau^{8/3}+3C_1\tau^2-3}}\ d\tau}+c_1,
  $$
  with $c_1\in \mathbb{R}$ a constant and $\xi_{00}\in \left(\xi_{01},\xi_{02}\right)$.
\end{itemize}
\end{theorem}

\begin{sproof}
The first item follows by standard arguments. For the second item, we note that choosing $\tilde{C}_1=3^{1/4}\cdot16C_1$ in \eqref{eq:Y_{C_1}} and using the change of coordinates $(\kappa,v)=\left(3^{-3/2}\xi^{4/3},\left(3^{-1/8}\sqrt{C_1}\theta\right)/4\right)$, the metric induced by $Y_{3^{1/4}\cdot 16C_1}$ coincides with $g_{C_1}$.

Then, we define $\phi_{C_1}$ as
$$
\phi_{C_1}(\xi,\theta)=Y_{3^{1/4}\cdot 16C_1}\left(3^{-3/2}\xi^{4/3},\frac{3^{-1/8}\sqrt{C_1}\theta}{4}\right).
$$
\end{sproof}

\begin{remark}
The limits $\lim_{\xi\searrow \xi_{01}} \zeta(\xi)=\zeta\left(\xi_{01}\right)$ and  $\lim_{\xi\nearrow \xi_{02}}\zeta(\xi)=\zeta\left(\xi_{02}\right)$ are finite.
\end{remark}

\begin{remark}
For simplicity, we choose $\xi_{00}=\left(9C_1/4\right)^{3/2}$.
\end{remark}

\begin{remark}
The immersion $\phi_{C_1}$ depends on the sign $\pm$ and on the constant $c_1$ in the expression of $\zeta$. As the classification is up to isometries of $\mathbb{S}^3$, the sign and the constant are not important, but they will play an important role in the gluing process.
\end{remark}


The construction of complete biconservative surfaces in $\mathbb{S}^3$ consists in two steps, and the key idea is to notice that $\left(D_{C_1},g_{C_1}\right)$ is, locally and intrinsically, isometric to a surface of revolution in $\mathbb{R}^3$.

The \textit{first step} is to construct a complete surface of revolution in $\mathbb{R}^3$ which on an open dense subset is locally isometric to $\left(D_{C_1},g_{C_1}\right)$. We start with the next result.

\begin{theorem}[\cite{N}]
\label{theorem3.18}
Let us consider $\left(D_{C_1},g_{C_1}\right)$ as above. Then $\left(D_{C_1},g_{C_1}\right)$ is the universal cover of the surface of revolution in $\mathbb{R}^3$ given by
\begin{equation}
\label{psi(xi,theta)}
\psi_{C_1,C_1^\ast}(\xi,\theta)=\left(\chi(\xi)\cos \frac{\theta}{C_1^\ast}, \chi(\xi)\sin \frac{\theta}{C_1^\ast},\nu(\xi)\right),
\end{equation}
where $\chi(\xi)=C_1^\ast/\xi$,
\begin{equation}
\label{h(xi)}
\nu(\xi)=\pm\bigintsss_{\xi_{00}}^\xi{\sqrt{\frac{3\tau^2-\left(C_1^\ast\right)^2\left(-\tau^{8/3}+3C_1\tau^2-3\right)}{\tau^4\left(-\tau^{8/3} +3C_1\tau^2-3\right)}}\ d\tau}+c_1^\ast,
\end{equation}
$C_1^\ast\in \left(0,\left(C_1-4/{3^{3/2}}\right)^{-1/2}\right)$ is a positive constant and $c_1^\ast\in\mathbb{R}$ is constant.
\end{theorem}

\begin{remark}
The immersion $\psi_{{C}_1, {C}_1^\ast}$ depends on the sign $\pm$ and on the constant $c_1^\ast$ in the expression of $\nu$. We denote by $S^{\pm}_{{C}_1,{C}_1^\ast,c_1^\ast}$ the image of $\psi_{{C}_1, {C}_1^\ast}$.
\end{remark}

\begin{remark}
The limits $\lim_{\xi\searrow \xi_{01}} \nu(\xi)=\nu\left(\xi_{01}\right)$ and  $\lim_{\xi\nearrow \xi_{02}}\nu(\xi)=\nu\left(\xi_{02}\right)$ are finite.
\end{remark}

We note that the boundary of $S^{\pm}_{{C}_1,{C}_1^\ast,c_1^\ast}$ is given by the curves
\begin{equation*}
\left(\frac{C_1^\ast}{\xi_{01}}\cos \frac{\theta}{C_1^\ast},\frac{C_1^\ast}{\xi_{01}}\sin \frac{\theta}{C_1^\ast}, \nu\left(\xi_{01}\right)\right)
\end{equation*}
and
\begin{equation*}
\left(\frac{C_1^\ast}{\xi_{02}}\cos \frac{\theta}{C_1^\ast},\frac{C_1^\ast}{\xi_{02}}\sin \frac{\theta}{C_1^\ast}, \nu\left(\xi_{02}\right)\right)
\end{equation*}
These curves are circles in affine planes in $\mathbb{R}^3$ parallel to the $Oxy$ plane and their radii are $C_1^\ast/\xi_{01}$ and $C_1^\ast/\xi_{02}$, respectively.

At a boundary point, using the coordinates $(\nu, \theta)$, we get that the tangent plane to the closure of $S^{\pm}_{{C}_1,{C}_1^\ast,c_1^\ast}$ is spanned by a vector which is tangent to the corresponding circle and by the vector $(0,0,1)$. Thus, the tangent plane is parallel to the rotational axis $Oz$.

Geometrically, we start with a piece of type $S^{\pm}_{{C}_1,{C}_1^\ast,c_1^\ast}$ and by symmetry to the planes where the boundary lie, we get our complete surface $\tilde{S}_{{C}_1,{C}_1^\ast}$; the process is periodic and we perform it along the whole $Oz$ axis.

Analytically, we fix $C_1$ and ${C}_1^\ast$, and alternating the sign and with appropriate choices of the constant $c_1^\ast$, we can construct a complete surface of revolution $\tilde{S}_{{C}_1,{C}_1^\ast}$ in $\mathbb{R}^3$ which on an open subset is locally isometric to $\left(D_{{C}_1},g_{{C}_1}\right)$. In fact, these choices of $+$ and $-$, and of the constants $c_1^\ast$ are uniquely determined by the ``first'' choice of $+$, or of $-$, and of the constant $c_1^\ast$. We start with $+$ and $c_1^\ast=0$.

The profile curve of $S^{\pm}_{{C}_1,{C}_1^\ast,c_1^\ast}$ can be seen as the graph of a function depending on $\nu$ and this allows us to obtain a function $F$ such that the profile curve of $\tilde{S}_{{C}_1,{C}_1^\ast}$ to be the graph of the function $\chi\circ F$ depending on $\nu$ and defined on the whole $Oz$ (or $O\nu$). The function $F:\mathbb{R}\to\left[\xi_{01},\xi_{02}\right]$ is periodic and at least of class $C^3$.

\begin{theorem}[\cite{N}]
\label{main_th3}
The surface of revolution given by
\begin{equation*}
\Psi_{C_1,C^\ast_1}(\nu,\theta)=\left((\chi\circ F)(\nu)\cos\frac{\theta}{C^\ast_1}, (\chi\circ F)(\nu)\sin \frac{\theta}{C^\ast_1},\nu \right), \qquad (\nu,\theta)\in\mathbb{R}^2,
\end{equation*}
is complete and, on an open dense subset, it is locally isometric to $\left(D_{C_1}, g_{C_1}\right)$. The induced metric is given by
$$
g_{C_1, C^\ast_1}(\nu,\theta)=\frac{3F^2(\nu)}{3F^2(\nu)- \left(C^\ast_1\right)^2(-F^{8/3}(\nu)+3C_1F^2(\nu)-3)}d\nu^2 +\frac{1}{F^2(\nu)}d\theta^2,
$$
$(\nu,\theta)\in\mathbb{R}^2$. Moreover, $\grad K\neq 0$ at any point of that open dense subset, and $1-K>0$ everywhere.
\end{theorem}

From Theorem $\ref{main_th3}$ we easily get the following result.

\begin{proposition}[\cite{N}]
The universal cover of the surface of revolution given by $\Psi_{C_1,C^\ast_1}$ is $\mathbb{R}^2$ endowed with the metric $g_{C_1,C^\ast_1}$. It is complete, $1-K>0$ on $\mathbb{R}^2$ and, on an open dense subset, it is locally isometric to $\left(D_{C_1},g_{C_1}\right)$ and $\grad K\neq 0$ at any point. Moreover any two surfaces $\left(\mathbb{R}^2,g_{C_1,C_1^\ast}\right)$ and $\left(\mathbb{R}^2,g_{C_1,C_1^{\ast\prime}}\right)$ are isometric.
\end{proposition}

The \textit{second step} is to construct effectively the biconservative immersion from $\left(\mathbb{R}^2,g_{C_1,C_1^\ast}\right)$ in $\mathbb{S}^3$, or from $\tilde{S}_{C_1,C_1^\ast}$ in $\mathbb{S}^3$. The geometric ideea of the construction is the following: from each piece $S^{\pm}_{{C}_1,{C}_1^\ast,c_1^\ast}$ of $\tilde{S}_{{C}_1,{C}_1^\ast}$ we ``go back'' to $\left(D_{C_1},g_{C_1}\right)$ and then, using $\phi_{C_1}$ and a specific choice of $+$ or $-$ and of the constant $c_1$, we get our biconservative immersion $\Phi_{C_1,C_1^\ast}$. Again, the choices of $+$ and $-$, and of the constant $c_1$ are uniquely determined  (modulo $2\pi$, for $c_1$) by the ``first'' choice of $+$, or of $-$, and of the constant $c_1$ (see \cite{N} for all details).

\vspace{0.5cm}

Some \textit{numerical experiments} suggest that $\Phi_{C_1,C_1^\ast}$ is not periodic and it has self-intersections along circles parallel to $Ox^3x^4$.

The projection of $\Phi_{C_1,C_1^\ast}$ on the $Ox^1x^2$ plane is a curve which lies in the annulus of radii $\sqrt{1-1/\left(C_1\xi_{01}^2\right)}$ and $\sqrt{1-1/\left(C_1\xi_{02}^2\right)}$. It has self-intersections and is dense in the annulus.

\vspace{0.5cm}

Concerning the biharmonic surfaces in $\mathbb{S}^3$ we have the following classification result.

\begin{theorem}[\cite{CMO01}]
Let $\varphi:M^2\to\mathbb{S}^3$ be a proper biharmonic surface. Then $\varphi(M)$ is an open part of the small hypersphere $\mathbb{S}^2(1/\sqrt{2})$.
\end{theorem}

\newpage
\appendix
\section{}
In the $c=0$ case, the idea was to construct, by symmetry, a complete biconservative surface in $\mathbb{R}^3$ starting with a piece of a biconservative surface. We illustrate this in the following figure obtained for $C_0=1$.

\begin{tikzpicture}
\node [inner sep=0pt] (local) at (0,0) {\includegraphics[scale=0.4]{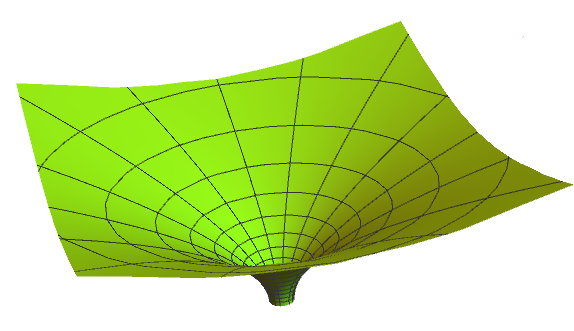}};
\node [inner sep=0pt] (global) at (10,0) {\includegraphics[scale=0.4]{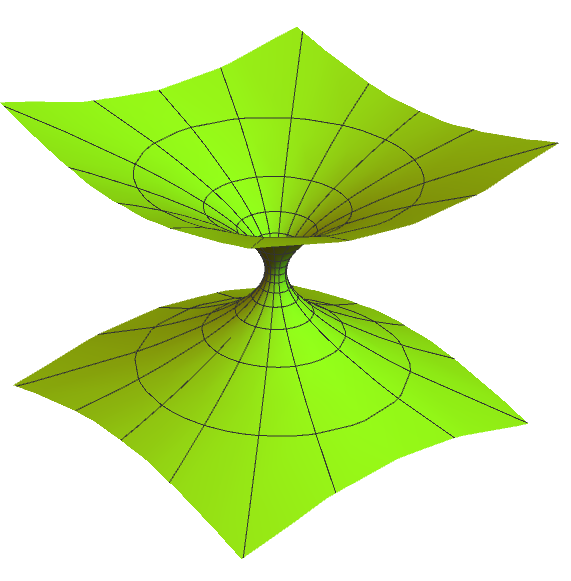}};
\draw[->,thin](local)--(global);
\end{tikzpicture}

\vspace{0.5cm}

In the $c=1$ case, the construction of a complete biconservative surface in $\mathbb{S}^3$ can be summarized in the next diagram, obtained for $C_1=C_1^\ast=1$, $c_1^\ast=0$ and we started with $+$ in the expression of $\nu$.

\begin{tikzpicture}
\node [inner sep=0pt] (manifold) at (3.6,0) {\includegraphics[scale=0.3]{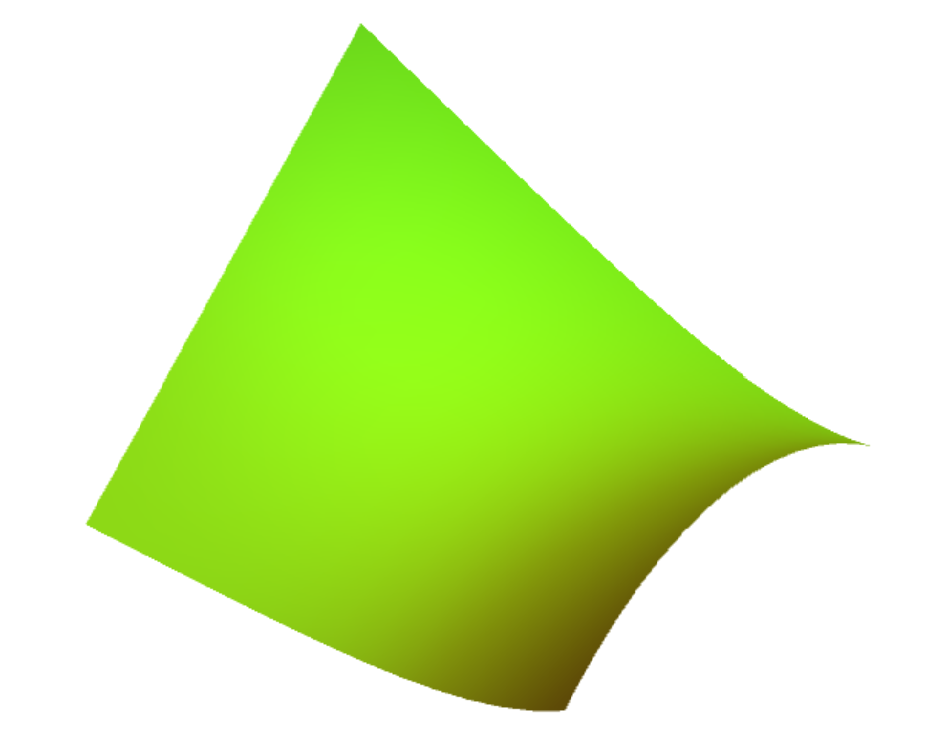}};
\draw (4.3,-2) node[below] {$\left(M^2,g\right)$};
\node [inner sep=0pt] (band) at (12,0) {\includegraphics[scale=0.4]{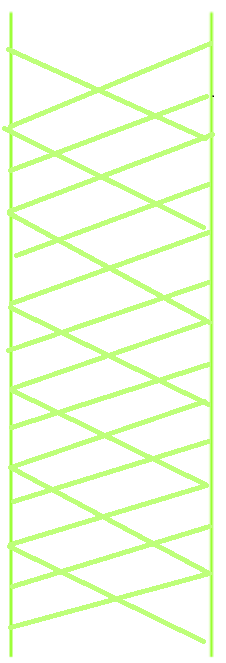}};
\draw [->,>=triangle 45] (10,0) -- (14,0);
\draw [->,>=triangle 45] (10.5,-2.5) -- (10.5,2.5);
\draw (11.2,0) node[below] {$\xi_{01}$};
\draw (12.9,0) node[below] {$\xi_{02}$};
\draw (14.1,0) node[below] {$\xi$};
\draw (10.5,2.6) node[left] {$\theta$};
\draw (12,-2.5) node[below] {$\left(D_{C_1},g_{C_1}\right)$};
\draw[<->,thick](manifold)-- node[midway,fill=none,above] {ISOMETRY}(9.5,0);
\draw[-,thick](11,-3.3)-| (11,-4);
\draw[-,thick](11,-4)-| (13,-4);
\draw[->,thick](13,-4) -- (13,-8) node[midway,sloped,fill=none,below] { $\phi_{C_1}=\phi^{\pm}_{C_1,c_1}$} node[midway,sloped,fill=none,above] {BICONSERVATIVE};
\node[ellipse, fill=brown, inner sep=0.2in] at (13,-9.3){$\mathbb{S}^3$};
\draw[<-,thick](10.7,-3.3)-| (10.7,-4);
\draw[->,thick](10.7,-4) --(3.5,-4) node[midway,sloped,fill=none,below] {$\psi_{C_1,C_1^\ast}=\psi^{\pm}_{C_1,C_1^\ast,c_1^\ast}$} node[midway,sloped,fill=none,above] {ISOMETRY} ;
\node[inner sep=0pt] at (2,-4){\includegraphics[scale=0.5]{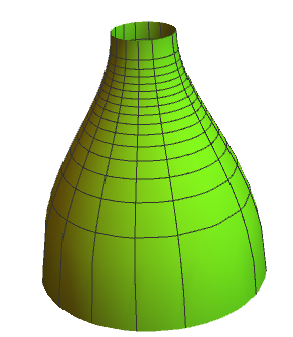}};
\draw (2,-5.5) node[below] {$S^{\pm}_{C_1,C_1^\ast,c_1^\ast}\subset\mathbb{R}^3$};
\node [inner sep=0pt] (graph6) at (5.3,-9.3) {\includegraphics[scale=0.6]{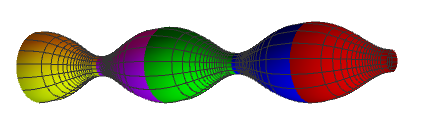}};
\draw (5.5,-10.2) node[below] {$\tilde{S}_{C_1,C_1^\ast}\subset\mathbb{R}^3$ \textcolor{red}{complete}};
\draw[->,thick](2,-6.2)--(5,-8)  node[pos=0.5,sloped, fill=none,above] {playing with the} node[pos=0.5,sloped,fill=none,below] {constant $c_1^\ast$ and $\pm$};
\begin{scope}[scale=2]
        \node (A) at (4,-4.2) {};
        \node (B) at (6.2,-4.1){};
        \draw[blue] plot [smooth,tension=1]
        coordinates {(A) (2.4,-3.2) (2.8,-2.7) (5.2,-2.6) (5.7,-3.6) (B)}
        [arrow inside={end=stealth,opt={red,scale=2}}{0.2,0.4,0.75,0.999}];
\end{scope}
\draw (6.5,-6.4) node[right] {$\begin{array}{c}\text{playing with the constant}\\ c_1 \text{ and } \pm \end{array}$};
\end{tikzpicture}

\vspace{0.5cm}

The projection of $\Phi_{1,1}$ on the $Ox^1x^2$ plane is represented in the next figure ($c_1=0$).

\begin{center}
\begin{tikzpicture}[xscale=1,yscale=1]

\node[inner sep=0pt]  at (.2,0){\includegraphics[width=.4\textwidth]{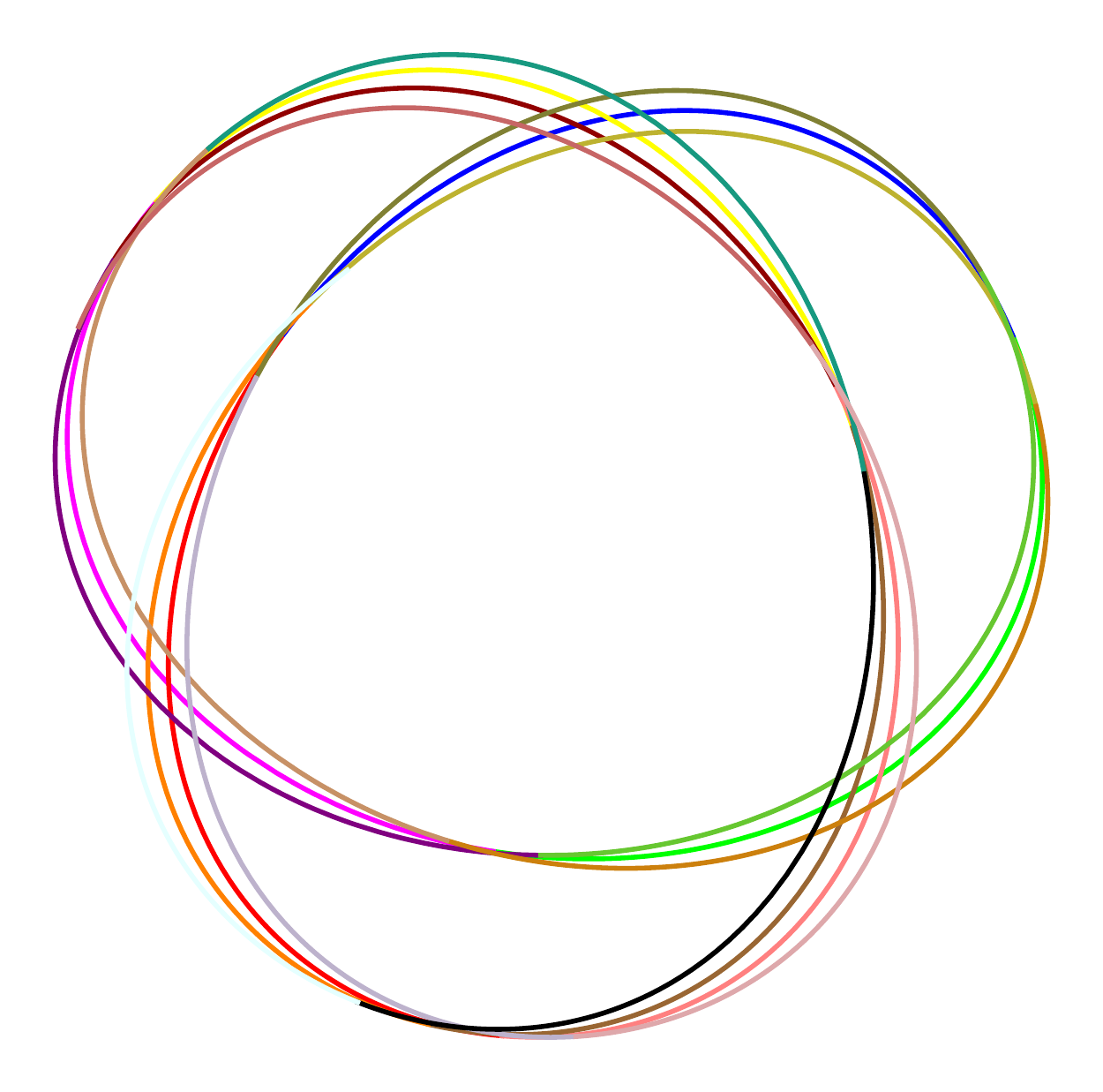}};
\draw [->,>=triangle 45] (-3,0) -- (3.5,0);
\draw [->,>=triangle 45] (0,-3) -- (0,3);
\draw (3.5,0) node[below] {$x^1$};
\draw (0,3) node[left] {$x^2$};
\end{tikzpicture}
\end{center}

\vspace{0.5cm}

The last two figures represent the signed curvature of the profile curve of $\tilde{S}_{C_1,C_1^\ast}$ and the signed curvature of the curve obtained projecting $\Phi_{1,1}$ on the $Ox^1x^2$ plane.

\begin{center}
\begin{tikzpicture}[xscale=1,yscale=1]

\node[inner sep=0pt]  at (.4,0){\includegraphics[width=.55\textwidth]{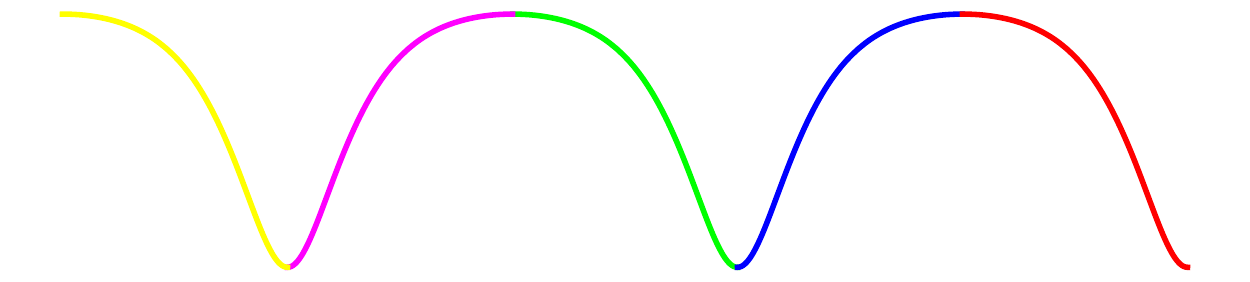}};
\draw [->,>=triangle 45] (-4.5,0.2) -- (5,0.2);
\draw [->,>=triangle 45] (0.85,-1.2) -- (0.85,1.5);
\draw (5,0.2) node[below] {$\nu$};
\draw (0.85,1.5) node[left] {$\kappa$};
\end{tikzpicture}
\end{center}

\begin{center}
\begin{tikzpicture}[xscale=1,yscale=1]

\node[inner sep=0pt]  at (.2,0.4){\includegraphics[width=.7\textwidth]{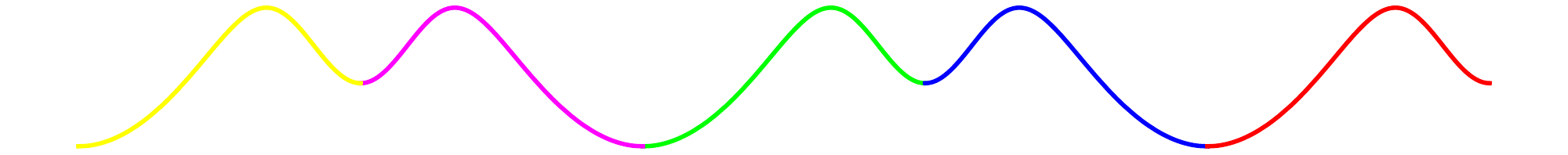}};
\draw [->,>=triangle 45] (-4.5,-0.5) -- (5,-0.5);
\draw [->,>=triangle 45] (0.7,-1) -- (0.7,2);
\draw (5,-0.5) node[below] {$\nu$};
\draw (0.7,2) node[left] {$\kappa$};
\end{tikzpicture}
\end{center}


\end{document}